# Basepoint freeness for nef and big line bundles in positive characteristic

By Seán Keel*

**Abstract**

A necessary and sufficient condition is given for semi-ampleness of a numerically effective (nef) and big line bundle in positive characteristic. One application is to the geometry of the universal stable curve over $\overline{M}_g$, specifically, the semi-ampleness of the relative dualizing sheaf, in positive characteristic. An example is given which shows this and the semi-ampleness criterion fail in characteristic zero. A second application is to Mori's program for minimal models of 3-folds in positive characteristic, namely, to the existence of birational extremal contractions.

### Introduction and statement of results

A map from a variety to projective space is determined by a line bundle and a collection of global sections with no common zeros. As all maps between projective varieties arise in this way, one commonly wonders whether a given line bundle is generated by global sections, or equivalently, if the associated linear system is basepoint free. Once a line bundle $L$ has a section, one expects the positive tensor powers $L^{\otimes n}$ to have more sections. If some such power is globally generated, one says that $L$ is *semi-ample*.

Semi-ampleness is particularly important in Mori's program for the classification of varieties (also known as the minimal model program, MMP). Indeed a number of the main results and conjectures — the Basepoint Free Theorem, the Abundance Conjecture, quasi-projectivity of moduli spaces — are explicitly issues of semi-ampleness. Some details will be given below.

There is a necessary numerical condition for semi-ampleness. The restriction of a semi-ample line bundle to a curve must have nonnegative degree; thus, if the line bundle $L$ on $X$ is semi-ample, then $L$ is *nef*; i.e., $L \cdot C \geq 0$ for every irreducible curve $C \subset X$. By a result of Kleiman (see [Kol96,VI.2.17]), nefness is equivalent to the apparently stronger condition: $L^{\dim Z} \cdot Z \geq 0$ for every proper irreducible $Z \subset X$. (We note this in relation to (0.1) below.)

---

*Partially supported by NSF grant DMS-9531940.

Nefness does not in general imply semi-ampleness (a nontorsion degree zero line bundle on a curve gives a counterexample).

The main result of this paper is a simple necessary and sufficient condition, in positive characteristic, for semi-ampleness of nef line bundles which are close to being ample. The statement involves a few natural notions:

0.0 *Definition-Lemma* (see [Kol96,VI.2.15,VI.2.16]). A line bundle $L$ on a scheme $X$, is called *big* if $L^{\otimes n}$ defines a birational *rational* map for $n >> 0$. If $L$ is nef, and $X$ is reduced and projective over a field, then $L$ is big if and only if $L^{\dim X_i} \cdot X_i > 0$ for any irreducible component $X_i$ of $X$. If $L$ is semi-ample, then $L$ is big if and only if the associated map is birational.

Associated to a nef and line bundle is a natural locus:

0.1 *Definition.* Let $L$ be a nef line bundle on a scheme $X$ proper over a field $k$. An irreducible subvariety $Z \subset X$ is called *exceptional* for $L$ if $L|_Z$ is not big, i.e. if $L^{\dim Z} \cdot Z = 0$. If $L$ is nef the *exceptional locus* of $L$, denoted by $\mathbb{E}(L)$, is the closure, with reduced structure, of the union of all exceptional subvarieties.

If $L$ is semi-ample, then $\mathbb{E}(L)$ is just the exceptional locus of the associated map. In (0.1) it is easy to check that one does not need to take the closure; $\mathbb{E}(L)$ is the union of finitely many exceptional subvarieties. (See (1.2).) Of course if $L$ is not big (and $X$ is irreducible), $\mathbb{E}(L) = X$. Observe that by Nakai's criterion for ampleness, [Kol96,VI.2.18], $L$ is ample if and only if $\mathbb{E}(L)$ is empty.

0.2 THEOREM. *Let $L$ be a* nef *line bundle on a scheme $X$, projective over a field of positive characteristic. $L$ is semi-ample if and only if $L|_{\mathbb{E}(L)}$ is semi-ample.*

In dimension two, (0.2) and the main ideas of its proof are contained in [Ar62,2-2.11]. (We received this unhappy news from Angelo Vistoli.) Similar questions were considered by Zariski; see [Z60].

0.3 COROLLARY (Assumptions as in (0.2)). *Assume the basefield is the algebraic closure of a finite field. If $L|_{\mathbb{E}(L)}$ is numerically trivial, then $L$ is semi-ample. In particular this holds if $\mathbb{E}(L)$ is one-dimensional. If $X$ is two-dimensional, any nef line bundle which is either big or numerically trivial is semi-ample.*

*Applications to* $\overline{\mathcal{M}}_{g,n}$. My first application of (0.2) is to the geometry of the universal stable pointed curve

$$\pi : \mathcal{U}_{g,n} \to \overline{\mathcal{M}}_{g,n}.$$

Let $\Sigma \subset \mathcal{U}_{g,n}$ be the union of the $n$ universal sections. Using (0.2) we prove:

0.4 THEOREM. $\omega_\pi(\Sigma)$ *is* nef *and big, and its exceptional locus is contained in the Deligne-Mumford boundary. If the basefield has positive characteristic, then $\omega_\pi(\Sigma)$ is semi-ample, but this fails in characteristic zero.*

The nefness and bigness of $\omega_\pi(\Sigma)$ were previously known; see, for example, [Spz81, p. 56], [V89], and [Kol90,4.6]. Note that even nefness is not obvious: the bundle is ample on fibres by the definition of stability, but *horizontal nefness* comes as a surprise. Bigness is related to additivity of the Kodaira dimension.

The failure in characteristic zero follows from a simple example, (3.0). The same example shows (0.2) fails in characteristic zero. It is also a counterexample to an interesting conjecture of Looijenga:

In Konsevich's proof of Witten's conjectures on the cohomology of $\overline{M}_g$ (see [Kon92]), a topological quotient $q : \overline{M}_{g,1} \to K_1 \mathcal{M}_g^1$ plays an important role. The topological space $K_1 \mathcal{M}_g^1$ is related to so-called Ribbon Graphs, and has other applications; see [HL96, §6]. In [L95], Looijenga raises the question of whether or not $K_1 \mathcal{M}_g^1$ admits any algebraic structure. He notes that the fibres of $q$ are algebraic, in fact (either points or) exactly the unions of exceptional curves for $\omega_\pi$. This implies that $\omega_\pi$ is semi-ample if and only if $K_1 \mathcal{M}_g^1$ is projective, and the quotient map $q : \overline{M}_{g,1} \to K_1 \mathcal{M}_g^1$ is the map associated to $\omega_\pi$. He conjectures that these statements hold (in characteristic zero). However example (3.0) implies that $q$ cannot even be a morphism of schemes. (See (3.6).) We do not know whether or not $K_1 \mathcal{M}_g^1$ (over $\mathbb{C}$) is an algebraic space (or equivalently, in terms of Definition (0.4.1) below, whether or not $\omega_\pi$ is Endowed With a Map (EWM).

*Applications to Mori's program.* The second application of (0.2) is to Mori's program for 3-folds, in positive characteristic.

We begin with some brief remarks on the general philosophy, so that the statements below will be intelligible. (For a detailed overview of the program, see [KMM87].) After we have stated our results, we will compare them with the existing literature.

In order to build moduli spaces of varieties one looks for a natural map to projective space, and thus for a natural semi-ample line bundle. On a general smooth variety, the only available line bundles are $\omega_X$ and its tensor powers. If $\omega_X$ is not nef, then, as noted above, one cannot hope for a map. Instead one looks for a birational modification which (morally speaking) increases the nefness of $\omega_X$. For surfaces one blows down a $-1$ curve. One focus of Mori's program is to generalize this procedure to higher dimensions. In this spirit, we have the next result, which is implied by the stronger but more technical results that follow:

COROLLARY (existence of extremal contraction). *Let $X$ be a projective $\mathbb{Q}$-factorial normal 3-fold defined over the algebraic closure of a finite field. Assume $X$ has nonnegative Kodaira dimension (i.e. $|mK_X|$ is nonempty for some $m > 0$). If $K_X$ is not nef, then there is a surjective birational map $f : X \to Y$ from $X$ to a normal projective variety $Y$, with the following properties*:

(1) $-K_X$ *is relatively ample.*
(2) *$f$ has relative Picard number one. More precisely: $f$ is not an isomorphism, any two fibral curves are numerically equivalent, up to scaling, and for any fibral curve $C$ the following sequence is exact*:

$$0 \longrightarrow \operatorname{Pic}(Y)_{\mathbb{Q}} \xrightarrow{f^*} \operatorname{Pic}(X)_{\mathbb{Q}} \xrightarrow{L \to L \cdot C} \mathbb{Q} \longrightarrow 0.$$

(Above, and throughout the paper, a *fibral curve* for a map, is a curve contained in a fibre.)

A more precise statement of these results involves a weakening of semi-ampleness:

0.4.1. *Definition.* A nef line bundle $L$ on a scheme $X$ proper over a field $k$ is *Endowed With a Map* (EWM) if there is a proper map $f : X \to Y$ to a proper algebraic space which contracts exactly the $L$-exceptional subvarieties, for example, for a subvariety $Z \subset X$, $\dim(f(Z)) < \dim(Z)$ if and only if $L^{\dim Z} \cdot Z = 0$. In particular, an irreducible curve $C \subset X$ is contracted by $f$ if and only if $L \cdot C = 0$. The Stein factorization of $f$ is unique; see (1.0).

For indications of the relationship between EWM and semi-ample, see (1.0).

We have a version of (0.2), (0.3) and the above corollary with EWM: whatever is stated for semi-ampleness over a finite field, holds for EWM over a field of positive characteristic. The exact statements are given in the body of the paper; see (1.9) and (1.9.1).

0.5 THEOREM (a basepoint-free theorem for big line bundles). *Let $X$ be a normal $\mathbb{Q}$-factorial three-fold, projective over a field of positive characteristic. Let $L$ be a* nef *and big line bundle on $X$.*

*If $L - (K_X + \Delta)$ is nef and big for some pure boundary $\Delta$ then $L$ is* EWM. *If the basefield is the algebraic closure of a finite field, then $L$ is semi-ample.*

(A boundary is a $\mathbb{Q}$-Weil divisor $\sum a_i D_i$ with $0 \le a_i \le 1$. It is called pure if $a_i < 1$ for all $i$.)

Note that when $K_X + \Delta$ has nonnegative Kodaira dimension, one does not need in (0.5) the assumption that $L$ is big, as it follows from the bigness of $L - (K_X + \Delta)$.

Basepoint-free theorems are related to the existence of extremal contractions, as follows: Suppose $K_X$ is not nef. Let $H$ be an ample divisor. Let $m$ be the infimum over rational numbers $\lambda$ such that $K_X + \lambda H$ is ample. $K_X + mH$ will be nef, and should be zero on some curve, $C$ (otherwise, at least morally, $K_X + mH$ would be ample and we could take a smaller $m$). $L := K_X + mH$ is semi-ample by the Basepoint Free Theorem. If $f : X \to Y$ is the associated map, then since $K_X + mH$ is pulled back from $Y$, $-K_X$ is relatively ample. There are, however, complications. For example, since we take the infimum, we need to show that $m$ is rational, so that some multiple of $K_X + mH$ is a line bundle. This leads to the study of the Mori-Kleiman Cone of Curves, $\overline{NE}_1(X)$, which is the closed convex cone inside $N_1(X)$ (the Neron-Severi group with real coefficients) generated by classes of irreducible curves. Specifically, one would like to know that the edges of the cone, at least in the half-space of $N_1(X)$ where $K_X$ is negative, are discrete, and generated by classes of curves. In this direction we have:

0.6 PROPOSITION (a cone theorem for $\kappa \geq 0$). *Let $X$ be a normal $\mathbb{Q}$-factorial three-fold, projective over a field. Let $\Delta$ be a boundary on $X$. If $K_X + \Delta$ has nonnegative Kodaira dimension, then there is a countable collection of curves $\{C_i\}$ such that*

(1) $\overline{NE}_1(X) = \overline{NE}_1(X) \cap (K_X + \Delta)_{\geq 0} + \sum_i \mathbb{R} \cdot [C_i].$

(2) *All but a finite number of the $C_i$ are rational and satisfy $0 < -(K_X + \Delta) \cdot C_i \leq 3$.*

(3) *The collection of rays $\{\mathbb{R} \cdot [C_i]\}$ does not accumulate in the half-space $(K_X)_{<0}$.*

In characteristic zero, (0.6) is contained in [Kol92,5.3].

The proof of (0.6) is simple; since we assume that (a multiple of) $K_X + \Delta$ is effective, the problem reduces to the cone theorem for surfaces.

*Brief overview of related literature.* For smooth $X$, over any basefield, Mori's original arguments, with extensions by Kollár, give much stronger results than mine. (See [Mo82], [Kol91].) The proofs are based on deformation theory, which (at least with current technology) requires very strong assumptions on the singularities. As one is mostly interested in smooth $X$, this may not at first seem like a serious restriction. However even if one starts with a smooth variety, the program may lead to singularities; for example, in the above corollary, $Y$ can have singularities, even if $X$ does not. Kollár has been able to extend the deformation methods to a fairly broad class of singularities, so-called LCIQ singularities. (See [Kol92].) In characteristic zero these include terminal 3-fold singularities, the singularities that occur in MMP be-ginning with smooth $X$, but this is not known in characteristic $p$. For other interesting applications of Kollár's main technical device, the bug-eyed cover, see [KMcK95].

In characteristic zero, with log terminal singularities, much stronger forms (without bigness assumptions) of all of the above results are known, in all dimensions. These are due to Kawamata and Shokurov; see [KMM87]. The proofs make essential use of vanishing theorems, which fail (at least in general) in positive characteristic.

In one important special case, that of a semi-stable family of surfaces, the full program is known in all (including mixed) characteristics. (See [Ka94].) In this case the cone and basepoint-free theorems are essentially surface questions, where the program is known in all characteristics; see [MT83]. Flips are another (very serious) matter.

My proof of (0.5) is based on ideas quite different from those of any these authors. It is a straightforward application of (0.2), and does not use any vanishing theorems or deformation theory. Note that we do not make any singularity assumptions of the log terminal sort.

*Overview of contents.* The proofs of (0.2) and (0.3) are in Section 1. Section 2 contains technical results about EWM and semi-ampleness used in the applications. A counterexample to (0.2), in characteristic zero, and various implications, are given in Section 3. In Section 4 we prove (0.4). The applications to Mori's program are in Section 5 which begins with the proof of (0.5). The proof of (0.6) is in (5.5).

*Thanks.* We received help on various aspects of this paper from a number of people. We would like to thank in particular S. Mori, J. McKernan, M. Boggi, L. Looijenga, F. Voloch, Y. Kawamata, R. Heitmann, K. Matsuki, M. Schupsky and A. Vistoli. We would like especially to thank J. Kollár for detailed and thoughtful (and occasionally sarcastic) comments on an earlier version of the paper.

0.9 *Notation and Conventions.* We will frequently mix the notation of line bundles and divisors. Thus, for example, if $L$ and $M$ are line bundles, we will write $L + M$ for $L \otimes M$.

We will often use the same symbol to denote a map, and a map induced from the map by applying a functor and we will sometimes denote the pullback $f^*(L)$ along a map $f : X \to Y$ by $L|_X$.

$X_{\text{red}}$ indicates the reduction of the space $X$. For a subspace $Y \subset X$, defined by an ideal sheaf $I \subset \mathcal{O}_X$, the $k^{\text{th}}$ *order neighborhood* is the subscheme defined by $I^{k+1}$.

For two Weil divisors $D, E$, we will say $D \geq E$ if the same inequality holds for every coefficient.

All spaces considered in this paper are assumed to be separated.

Whenever we have a basefield $k$, we implicitly assume maps between $k$-spaces are $k$-linear. The only non $k$-linear map which we consider in the paper is the ordinary Frobenius, defined below.

*Frobenius Maps.* For a scheme, $X$, of characteristic $p > 0$, and $q = p^r$, we indicate by $F_q : X \to X$ the (ordinary) Frobenius morphism, which is given by the identity on topological spaces, and the $q^{\text{th}}$ power on functions. (See [H77, IV.2.4.1].) If $X$ is defined over $k$, $F_q$ factors as

$$F_q : X \to X^{(q)} \to X$$

where $X^{(q)} \to X$ is the pullback of $F_q$ on $\mathrm{Spec}(k)$. The map $X \to X^{(q)}$ is called the geometric Frobenius. It is $k$-linear. When $X$ is of finite type over $k$, the Frobenius, and geometric Frobenius, are finite universal homeomorphisms; see [Kol95,§6].

## 1. Proof of main theorem

We begin with some simple properties of the map associated to an EWM line bundle.

1.0 DEFINITION-LEMMA. *Let $L$ be a* nef EWM *line bundle on an algebraic space $X$, proper over a field. There is a proper map $X \to Y$ as in (0.4.1.) (i.e. a map which contracts exactly the $L$-exceptional irreducible subspaces) related to $L$. If $f$ is such a map, and $f_*(\mathcal{O}_X) = \mathcal{O}_Y$, then $f$ is unique. We will call it the* map associated to $L$. *The associated map, $f$, has the following properties*:

(1) *If $f' : X \to Y'$ is a proper map, which contracts any proper irreducible curve $C \subset X$ with $L \cdot C = 0$, then $f'$ factors uniquely through $f$.*
(2) *$f'$ as in (1). If, in addition, $L$ is $f'$ numerically trivial then the induced map $Y \to Y'$ is finite, and $f$ is the Stein factorization of $f'$.*
(3) *$f$ is the Stein factorization of any map related to $L$.*
(4) *Let $h : X' \to X$ be any proper map. The pullback $L|_{X'}$ is EWM, and the associated map is the Stein factorization of $f \circ h$.*
(5) *$L$ is semi-ample if and only if $L^{\otimes m}$ is the pullback of a line bundle on $Y$ for some $m > 0$.*

*Proof.* Of course uniqueness follows from the universal property, (1), which in turn follows from the rigidity lemma, [Kol96,II.5.3]. The remaining remarks are either easy, or contained in (1.1) and (1.3) below. □

1.1 LEMMA. *Let $f: X \to Y$ be a proper surjective map with geometrically connected fibres, between algebraic spaces of finite type over a field. Assume either that the characteristic is positive, or that $f_*(\mathcal{O}_X) = \mathcal{O}_Y$. If $L$ is semi-ample and $f$-numerically trivial, then, for some $r > 0$, $L^{\otimes r}$ is pulled back from a line bundle on $Y$.*

*Proof.* Let $f': X \to Y'$ be the Stein factorization of $f$. Then by assumption, $Y' \to Y$ is a finite universal homeomorphism (the identity in characteristic zero). By (1.4) it is enough to show $L^{\otimes r}$ is pulled back from $Y'$. Let $g: X \to Z$ be the map associated to $L$. Since $L$ is pulled back from an ample line bundle on $Z$, every fibre of $f'$ is contracted by $g$. Thus $f'$ factors through $g$ by the rigidity lemma, [Kol96,II.5.3]. □

1.2. By [DG, 4.3.4], for proper map $f: X \to Y$, if $f_*(\mathcal{O}_X) = \mathcal{O}_Y$, then $f$ has geometrically connected fibers.

1.3 COROLLARY. *Let $L$ be a nef line bundle on an algebraic space $X$, proper over a field. Assume $L$ is EWM and $f: X \to Y$ is the associated map. The following are equivalent:*

(1) *$L$ is semi-ample.*

(2) *$L^{\otimes r}$ is pulled back from a line bundle on $Y$, for some $r > 0$.*

(3) *$L^{\otimes r}$ is pulled back from an ample line bundle on $Y$, for some $r > 0$.*

*Proof.* (1) implies (2) by (1.1) and (3) obviously implies (1). Thus it is enough to show (2) implies (3). Assume $L = f^*(M)$ for a line bundle $M$ on $Y$. Let $W \subset Y$ be an irreducible subspace of dimension $k$. Since $f$ is surjective, there is an irreducible $k$-dimensional subspace $W' \subset X$ surjecting onto $W$. Let $d$ be the degree of $W' \to W$. Since $W'$ is not $L$-exceptional,

$$0 < L^k \cdot W' = d(M^k \cdot W).$$

So $M$ is ample by Nakai's criterion, [H70]. □

The next two lemmas point up the advantage of positive characteristic.

1.4 LEMMA. *Let $f: X \to Y$ be a finite universal homeomorphism between algebraic spaces of finite type over a field $k$ of characteristic $p > 0$. Then for some $q = p^r$ the following hold, with $L$ a line bundle on $Y$:*

(1) *For any section $\sigma$ of $f^*(L)$, $\sigma^{\otimes q}$ is in the image of*

$$f^*: H^0(Y, L^{\otimes q}) \to H^0(X, f^*(L^{\otimes q})).$$

(2) *$L$ is semi-ample if and only if $f^*(L)$ is semi-ample.*

(3) *The map*
$$f^* : \operatorname{Pic}(Y) \otimes \mathbb{Z}[1/q] \to \operatorname{Pic}(X) \otimes \mathbb{Z}[1/q]$$
*is an isomorphism.*

(4) *If $f^*(\sigma_1) = f^*(\sigma_2)$, for two sections $\sigma_i \in H^0(Y, L)$, then $\sigma_1^{\otimes q} = \sigma_2^{\otimes q}$.*

*Proof.* By [Kol95,6.6] there are a finite universal homeomorphism $g : Y \to X$, and a $q$, as in the statement of the lemma, such that the composition $g \circ f$ is the Frobenius morphism, $F_q$; see (0.9). The map induced by $F_q$ on Cartier divisors is just the $q^{\text{th}}$ power. The result follows easily. □

There is also a version of (1.4) for EWM:

1.5 LEMMA. *Let $g : X \to X'$ be a finite universal homeomorphism between algebraic spaces proper over a field of positive characteristic. A line bundle $L$ on $X'$ is EWM if and only if $g^*(L)$ is EWM.*

*Proof.* Let $f : X \to Z$ be the map associated to $g^*(L)$. By [Kol95,6.6] there is a pushout diagram

$$\begin{array}{ccc} X & \xrightarrow{g} & X' \\ f \downarrow & & \downarrow f' \\ Z & \xrightarrow{\tilde{g}} & Z' \end{array}$$

with $\tilde{g}$ a finite universal homeomorphism. Clearly $f'$ is a related map for $L$. □

The main step in proving (0.2) is the following:

1.6 PROPOSITION. *Let $X$ be a projective scheme over a field of positive characteristic. Let $L$ be a nef line bundle on $X$. Suppose $L = A + D$ where $A$ is ample and $D$ is effective and Cartier. $L$ is EWM if and only if $L|_{D_{\text{red}}}$ is EWM. $L$ is semi-ample if and only if $L|_{D_{\text{red}}}$ is semi-ample.*

*Proof.* Assume $L|_{D_{\text{red}}}$ is EWM.

Let $D_k$ be the $k^{\text{th}}$ order neighborhood of $D$. By (1.5), $L|_D$ is EWM. Let $p : D \to Z$ be the associated map. Let $I = I_D$. Note that since $L|_D$ is numerically $p$-trivial,
$$D|_D = L|_D - A|_D$$
is $p$-anti-ample. Thus by Serre vanishing (and standard exact sequences) there exists $n > 0$ such that:

I. $R^i p_*(I^j/I^t) = 0$ for any $t \geq j \geq n$, $i > 0$.

II. Let $J = I^s$ for some $s > 0$. For any coherent sheaf $\mathcal{F}$ on $X$, $R^i p_*(J^k \cdot \mathcal{F}/J^{k+1} \cdot \mathcal{F})$ vanishes for $k >> 0$, $i > 0$.

III. For $J = I^n$, $p_*(\mathcal{O}/J^t) \to p_*(\mathcal{O}/J)$ is surjective for $t \geq 1$.

By (1.5), $L|_{D_n}$ is EWM. Let $D_n \to Z_n$ be the associated map. By (1.0) the induced map $p' : D \to Z_n$ factors through $p$ and the induced map $Z \to Z_n$ is finite. Thus $p'$ satisfies I–III. Replace $p$ by $p'$, and $Z$ by the scheme-theoretic image of $p'$. We will make similar adjustments to $p$ later in the argument, without further remark.

By II–III and [Ar70,3.1,6.3], there is an embedding $Z_n \subset X'$ of $Z_n$ in a proper algebraic space, and an extension of $D_n \to Z_n$ to a proper map $p : X \to X'$, such that $D$ is set-theoretically the inverse image of $Z$, and such that $p : X \setminus D \to X' \setminus Z$ is an isomorphism. Passing to the Stein factorization, we may assume $p_*(\mathcal{O}_X) = \mathcal{O}_{X'}$. From (1.7) it follows that $L$ is EWM and $p : X \to X'$ is the associated map.

Now suppose $L|_{D_{\text{red}}}$ is semi-ample. By (1.1) $L^{\otimes r(k)}|_{D_k}$ is pulled back from the scheme-theoretic image of $D_k$ in $X'$ for some $r(k) > 0$. Thus for some $r > 0$, $L^{\otimes r}$ is pulled back from $X'$ by (I) and (1.10). Thus $L$ is semi-ample by (1.3). □

1.7 LEMMA. *Let $L$ be a* nef *line bundle on a scheme $X$ proper over a field. If $L = A + D$ with $A$ ample and $D$ effective and Cartier, then $\mathbb{E}(L) \subset D$. In any case $\mathbb{E}(L)$ is a finite union of exceptional subvarieties.*

*Proof.* For the first claim, let $Z \subset X$ be an irreducible subvariety of dimension $k$. If $Z \not\subset D$, then $D|_Z$ is effective and Cartier. $L|_Z = A|_Z + D|_Z$. Thus
$$L^k \cdot Z \geq A^k \cdot Z > 0.$$

We will prove the second claim by induction on the dimension of $X$, and we can assume $L$ is big. We may write $L = A + D$ as in the statement, by Kodaira's lemma, [Kol96,VI.2.16]. Write $D_{\text{red}} = D_B + D_E$ where $D_B$ is the union of the irreducible components on which $L$ is big, and $D_E$ is the union of the remaining components. By the first claim $\mathbb{E}(L) = D_E \cup \mathbb{E}(L|_{D_B})$. □

1.8 LEMMA. *Let $X$ be an algebraic space, proper over a field of positive characteristic, and let $L$ be a* nef *line bundle on $X$. Suppose $X$ is union of closed subspaces $X = X_1 \cup X_2$. If $L|_{X_i}$ is semi-ample (resp. EWM) for $i = 1, 2$ and $\mathbb{E}(L) \subset X_1$ then $L$ is semi-ample (resp. EWM).*

*Sketch of Proof.* This follows easily from [Kol95,8.4] and [Ar70,6.1]. We will state and prove more general results in Section 2, which the reader may consult for a complete proof. See, for example, (2.10.1) and (2.12). □

1.9 THEOREM. *Let $L$ be a* nef *line bundle on a scheme $X$, projective over a field of positive characteristic. $L$ is semi-ample (resp. EWM) if and only if $L|_{\mathbb{E}(L)}$ is semi-ample (resp. EWM).*

*Proof.* We induct on the dimension of $X$. By (1.8), we may assume $L$ is big, and, by (1.4) and (1.5), we may assume $X$ is reduced. By Kodaira's lemma, [Kol96,VI.2.16], $L = A + D$ as in (1.6). $L|_D$ is semi-ample (resp. EWM) by induction. Now apply (1.6). □

1.9.1 COROLLARY. *With $(L, X)$ as in* (1.9), *if $L|_{\mathbb{E}(L)}$ is numerically trivial* (*in particular if $\mathbb{E}(L)$ is one-dimensional*) *then $L$ is* EWM, *and semi-ample if the basefield is the algebraic closure of a finite field.*

*Proof.* Note that any numerically trivial line bundle is EWM. Now, (1.9.1) follows from (1.9) and (2.16). □

1.10 LEMMA. *Let $f : X \to X'$ be a proper map between algebraic spaces, with $f_*(\mathcal{O}_X) = \mathcal{O}_{X'}$. Let $D \subset X$ be a subspace with ideal sheaf $I \subset \mathcal{O}_X$ and scheme-theoretic image $Z \subset X'$. Let $D_k \subset X$ be the $k^{\text{th}}$ order neighborhood of $D$. Assume $D$ is set-theoretically the inverse image of $Z$, and that the map*

$$f : X \setminus D \to X' \setminus Z$$

*is an isomorphism. Let $L$ be a line bundle on $X$, such that for each $k$ there is an $r(k) > 0$ such that $L^{\otimes r(k)}|_{D_k}$ is pulled back from (the scheme-theoretic image) $f(D_k) \subset X'$.*

*If $R^1 f_*(I^k/I^{k+1}) = 0$ for $k >> 0$, then $L^{\otimes r}$ is pulled back from $X'$ for some $r \geq 1$.*

*Proof.* Replace $L$ by a power so that $L|_D$ is pulled back from $Z$. Choose $n$ so that

(1.10.1) $$R^1 f_*(I^k/I^{k+1}) = 0$$

for $k \geq n$. Let $Z_i \subset X'$ be the scheme-theoretic image of $D_i$. Replace $L$ by $L^{\otimes r}$, so that $L|_{D_n}$ is pulled back from $Z_n$. We will show that:

(1) $f_*(L)$ is locally free of rank one, and
(2) the canonical map $f^*(f_*(L)) \to L$ is an isomorphism.

(1) and (2) can be checked after a faithfully flat extension of $X'$. The vanishing (1.10.1), and the assumption $f_*(\mathcal{O}_X) = \mathcal{O}_{X'}$ are preserved by such an extension. (1) and (2) are local questions along $Z$. Thus we may assume $X'$ is the spectrum of a local ring, $A$. It follows that $L|_{D_n}$ is trivial. By (1.6.1), we can choose, for all $i \geq 1$, global sections $\sigma_i \in H^0(L \otimes D_i)$ such that $\sigma_i|_{D_j} = \sigma_j$ for $i \geq j$, and such that $\sigma_1$ is nowhere vanishing. By Nakayama's lemma, $\sigma_i$ is nowhere vanishing and $L|_{D_i}$ is trivial, for all $i$. The collection $\{\sigma_i\}$ induces an isomorphism

$$\varprojlim H^0(L \otimes \mathcal{O}_{D_i}) \to \varprojlim H^0(\mathcal{O}_{D_i}).$$

By the theorem on formal functions, [DG], the left-hand side is $H^0(L) \otimes \hat{A}$ and the right-hand side is $H^0(\mathcal{O}_X) \otimes \hat{A} = \hat{A}$, where $\hat{A}$ is the completion of $A$ along $Z$. Hence we have (1).

By (1), to establish (2) we need only show surjectivity, or equivalently (in the current local situation), that $L$ is basepoint free. For any $k \geq n$

$$H^1(L \otimes I^k) \otimes \hat{A} = \varprojlim_r H^1(L \otimes I^k/I^{k+r}) \text{ by the theorem on formal functions}$$

$$= \varprojlim_r H^1(I^k/I^{k+r}) \text{ since } L|_{D_i} \text{ is trivial for all } i$$

$$= 0 \text{ by (1.10.1) and induction.}$$

Thus $H^0(L) \to H^0(L \otimes \mathcal{O}_{D_n})$ is surjective. Since $L|_{D_n}$ is trivial, (2) follows. □

1.10.2 *Remark.* The proof shows that if $R^1 f_*(I^k/I^{k+1}) = 0$ for all $k \geq n$, then $L^{r(n)}$ is a pull back. In characteristic $p$, by (1.4), one need only assume $L|_D$ is pulled back, and in this case $r(n)$ can be chosen independent of $L$. This strengthening of (1.10) was pointed out to me by the referee.

## 2. Descending EWM or semi-ample in pushout diagrams

Suppose there is a proper map $h : Y \to X$, and that $L$ is a nef line bundle on $X$. If $L$ is EWM, or semi-ample, then it easily follows (see (1.0)) that the same holds for $h^*(L)$. An important technical problem in the proofs of the main results of the paper will be to find conditions under which the reverse implication holds, that is, conditions under which EWM, or semi-ampleness, descends from $h^*(L)$ to $L$. A simple example is (1.8). This section contains a number of results of this sort.

2.1 LEMMA. *Let $i : Z \to X$ be a proper, set-theoretic surjection between algebraic spaces of finite type over a perfect field $k$ of positive characteristic. Let $f, g : X \to Y$ be maps such that $f \circ i = g \circ i$. Then there is a finite universal homeomorphism (over $k$) $h : Y \to Y'$ such that $h \circ f = h \circ g$.*

*Proof.* Replacing $Z$ by its scheme-theoretic image, we can assume $i$ is a closed embedding, defined by a nilpotent ideal. For some $q = p^r$, $F_q$ (on $X$) factors through $i$. By the functorality of $F_q$, it follows that $F_q \circ f = F_q \circ g$ (here $F_q$ is on $Y$).

$F_q$ factors as $Y \to Y^{(q)} \to Y$, where the first map is the ($k$-linear) geometric Frobenius, and the second is the pullback of $F_q$ on $\mathrm{Spec}(k)$. When the basefield is perfect, $Y^{(q)} \to Y$ is an isomorphism. Thus for $Y \to Y'$ we can take the geometric Frobenius. □

2.2 LEMMA. *Let $L$ be a line bundle on an algebraic space $X$ proper over a field $k$. If $k \subset k'$ is an algebraic field extension, then $L|_{X_{k'}}$ is semi-ample (resp. EWM) if and only if $L$ is semi-ample (resp. EWM), where $X_{k'} := X \times_k k'$.*

*Proof.* The semi-ample case is obvious. The EWM case follows from flat descent; see [Ar68,7.2]. □

2.3 *Notation.* We assume the following commutative diagram:

(2.4)
$$\begin{array}{ccc} \mathcal{C} & \xrightarrow{j} & Y \\ p \downarrow & & p \downarrow \\ \mathcal{D} & \xrightarrow{i} & X \end{array}$$

which is a pushout; i.e., for any scheme $T$, the induced diagram

(2.5)
$$\begin{array}{ccc} \operatorname{Hom}(X,T) & \xrightarrow{\circ i} & \operatorname{Hom}(\mathcal{D},T) \\ \circ p \downarrow & & \circ p \downarrow \\ \operatorname{Hom}(Y,T) & \xrightarrow{\circ j} & \operatorname{Hom}(\mathcal{C},T) \end{array}$$

is a pullback diagram (of sets). We assume $i$ and $j$ are closed embeddings, $p$ is proper and that for any open subset $U \subset X$, the diagram induced from (2.4) by restriction to $U$ is again a pushout. All the spaces are algebraic spaces of finite type over a field of positive characteristic.

Let $L$ be a nef line bundle on $X$. For any proper map $T \to X$ such that $L|_T$ is EWM, we denote the associated map by $g_T : T \to Z_T$, occasionally dropping subscripts when they are clear from context.

2.6 LEMMA (notation as in (2.3)). *Assume the basefield is perfect. Let $D \subset X$ be a reduced subspace, and $C \subset Y$ the reduction of its inverse image. Assume $\mathcal{D}$ is contained set-theoretically in $D$ and that $L|_D$ and $L_Y$ are EWM, and that $g_Y|_C$ has geometrically connected fibres. Then $L$ is EWM. If, furthermore, $L|_D$ and $L|_Y$ are semi-ample, then $L$ is semi-ample.*

*Proof.* We will replace $L$ at various times by a positive tensor power, often without remark.

By (1.0), $g_C : C \to Z_C$ is the Stein factorization of either

$$g_Y|_C : C \hookrightarrow Y \xrightarrow{g_Y} Z_Y$$

or

$$g_D \circ p|_C : C \xrightarrow{p} D \xrightarrow{g_D} Z_D.$$

Let $V := g_Y(C) \subset Z_Y$ be the scheme-theoretic image of the first map. Since $g_Y|_C$ has geometrically connected fibres, the induced map $Z_C \to V$ is a finite universal homeomorphism. Let $\overline{C}, \overline{D}$ be the reductions of $\mathcal{C}$ and $\mathcal{D}$.

Consider first the semi-ample case. We will use the following notation. If $T \to X$ is a map and $L|_T$ is semi-ample, then (after replacing $L$ by a power), $L|_T$ is pulled back from an ample line bundle on $Z_T$, by (1.0.5). We will write $L|_T = g_T^*(M_T)$.

The induced diagram (2.5), with $T = \mathbb{A}^1$, gives a short exact sequence of sheaves
$$0 \to \mathcal{O}_X \to p_*(\mathcal{O}_Y) \oplus \mathcal{O}_\mathcal{D} \xrightarrow{j^* \oplus -p^*} p_*(\mathcal{O}_\mathcal{C}).$$
There is an analogous exact sequence for global sections of a line bundle

(2.7) $\quad 0 \to H^0(X, L) \to H^0(Y, L|_Y) \oplus H^0(\mathcal{D}, L|_\mathcal{D}) \xrightarrow{j^* \oplus -p^*} H^0(\mathcal{C}, L|_\mathcal{C}).$

Choose a point $x \in X$. Observe that
$$d := g_D(i^{-1}(x)) \subset Z_D, \text{ (which is } \emptyset \text{ if } x \notin D) \text{ and}$$
$$y := g_Y(p^{-1}(x)) \subset Z_Y$$
are zero-dimensional. As $M_D$ is ample, there is a section $\sigma^{Z_D} \in H^0(Z_D, M_D)$, nonvanishing at any point of $d$. Let $\sigma^D, \sigma^C$ be the pullbacks of $\sigma^{Z_D}$ to $D, C$. By (1.4), replacing the sections by powers, we can assume $\sigma^C$ is pulled back from a section $\sigma^V \in H^0(V, M_Y|_V)$. Since $M_Y$ is ample, again replacing sections by powers, we may assume $\sigma^V$ is the restriction of $\sigma^{Z_Y}$, nonvanishing at any point of $y$. Let $\sigma^Y$ be the pullback of $\sigma^{Z_Y}$. Let $\sigma^{\overline{D}}$ be the restrictions of $\sigma^D$ to $\overline{D}$. By (1.4.1), replacing sections by powers, we can assume $\sigma^{\overline{D}}$ is the restriction of a section $\sigma^\mathcal{D}$ in $H^0(\mathcal{D}, L|_\mathcal{D})$. By construction $\sigma^Y$ and $p^*(\sigma^\mathcal{D})$ have the same restrictions to $\overline{C}$. Thus by (1.4.4), after replacing sections by powers, we can assume
$$\sigma^Y|_\mathcal{C} = p^*(\sigma^\mathcal{D})$$
and thus by (2.7), the sections join to give a global section $\sigma \in H^0(X, L)$. By construction, $\sigma$ does not vanish at $x$. $\square$

Now for the EWM part. Since the induced map $Z_C \to V$ is a finite universal homeomorphism, by [Kol95,8.4] there is a pushout diagram

$$\begin{array}{ccc} Z_C & \longrightarrow & V \\ \downarrow & & \downarrow \\ Z_D & \longrightarrow & Z_1 \end{array}$$

where the base is a finite universal homeomorphism. Since the left-hand column is finite, so is the right-hand column. By [Ar70,6.1] there is another pushout diagram

$$\begin{array}{ccc} V & \longrightarrow & Z_Y \\ \downarrow & & \downarrow \\ Z_1 & \longrightarrow & Z_2 \end{array}$$

where the base map is a closed embedding, and the right column is finite. Define $f, h$ to be the compositions
$$f : Y \xrightarrow{g_Y} Z_Y \to Z_2$$
$$h : D \xrightarrow{g_D} Z_D \to Z_1 \subset Z_2.$$

Note, by construction, that these are related maps for $L|_Y, L|_D$, and $h \circ p|_{\overline{C}} = f|_{\overline{C}}$. By [Kol95,8.4], there is a finite universal homeomorphism $k : Z_2 \to Z_3$, such that $k \circ h|_{\overline{D}}$ extends to a map $h' : \mathcal{D} \to Z_3$. Also, $h' \circ p$ and $k \circ f|_{\mathcal{C}}$ agree on $\overline{C}$; thus by (2.1), after replacing $k$ by its composition with another finite universal homeomorphism, we may assume
$$p \circ h' = k \circ f \circ j \text{ (on } \mathcal{C}\text{)}.$$
Now by (2.5), the maps join to give a map $g : X \to Z_3$ and $g \circ p$ agrees with $g_Y$ up to composition by a finite map, and hence $g \circ p$ is a related map for for $L|_Y$. Thus by (1.14), $g$ is a related map for $L$, and $L$ is EWM.  □

2.8 LEMMA. *Let $p : Y \to X$ be a proper surjection between algebraic spaces, proper over a field $k$. Let $g : X \to Z$ be a proper map. $L|_Y$ is EWM and $g \circ p$ is a related map, if and only if $L$ is EWM and $g$ is a related map.*

*Proof.* This is an easy consequence of the projection formula; see, for example, the proof of (1.3).  □

2.9 COROLLARY. *Let $X$ be a reduced algebraic space, proper over a field of positive characteristic. Suppose $X$ is union of closed subspaces $X = X_1 \cup X_2$. Let $L$ be a nef line bundle on $X$ such that $L|_{X_i}$ is EWM for $i = 1, 2$. Let $g$ be the map associated to $L|_{X_2}$. Assume $g|_{X_1 \cap X_2}$ has connected geometric fibres. Then $L$ is EWM. If furthermore $L|_{X_i}$ is semi-ample, for $i = 1, 2$, then $L$ is semi-ample.*

*Proof.* By (2.2) we can assume the basefield is algebraically closed. The diagram
$$\begin{array}{ccc} X_1 \cap X_2 & \longrightarrow & X_2 \\ \downarrow & & \downarrow \\ X_1 & \longrightarrow & X \end{array}$$
is a pushout diagram (where $X_1 \cap X_2$ is the scheme-theoretic intersection). So the result follows from (2.6).  □

2.10 LEMMA. *Let $p : Y \to X$ be a proper surjection between reduced algebraic spaces of finite type over a field of positive characteristic. Let $D \subset X$ be a reduced subspace, and $C \subset Y$ the reduction of its inverse image. Let $L$ be a line bundle on $X$ such that $L|_D$ and $p^*(L)$ are semi-ample. Let $g$ be the map associated to $p^*(L)$. Assume $g|_C$ has geometrically connected fibres.*

*If $X$ is normal outside of $D$, then $L$ is semi-ample.*

2.10.1 *Remark.* Note that the connectivity assumption of (2.10) is trivially satisfied if $p^*(L)$ is nef and big, and its exceptional locus is contained in $C$, for in that case, any fibre of $g$ is either a single point, or contained in $C$.

*Proof of* (2.10). Let $x \in X$ be a point. We note first that the proof of the semi-ample part of (2.6) yields (after replacing $L$ by a power) sections

$$\sigma^D \in H^0(D, L|_D)$$
$$\sigma^Y \in H^0(Y, p^*(L)), \text{ not vanishing at any point of } p^{-1}(x),$$

which have the same pullbacks to $\mathcal{C}$.

We will construct a section of (some power of) $L$ nonvanishing at $x$ and will use only the existence of the above sections, and the normality of $X \setminus D$ (but no assumptions of connectivity).

Let $\tilde{X} \to X$ be the normalization of $X$, with conductors $\mathcal{C} \subset \tilde{X}$, $\mathcal{D} \subset X$. We may assume $D$ is the reduction of $\mathcal{D}$.

Replacing $Y$ by a subspace, normalizing, and then taking a normal closure (in the sense of Galois theory), and replacing $\mathcal{C}$ and the sections by their pullbacks, we may assume $Y$ is normal, and that $p$ is generically finite and factors as

$$Y \xrightarrow{g} X' \xrightarrow{j} \tilde{X} \longrightarrow X$$

where $X'$ is normal, $g$ is the quotient by a finite group, and $j$ is a finite universal homeomorphism. Indeed, once $[K(Y) : K(X)]$ is a normal field extension, with Galois group $G$, define $X'$ to be the integral closure of $X$ in $K(Y)^G$.

We may replace $\sigma^Y$ by the tensor power of its translates under $G$ (which replaces the restriction to $\mathcal{C}$ by some tensor power), and so may assume $\sigma^Y$ is $G$ invariant. Then it is pulled back from $X'$. After replacing sections by powers, by (1.4), we may assume $Y = \tilde{X}$.

By [R94,2.1], the diagram

(2.11)
$$\begin{array}{ccc} \mathcal{C} & \longrightarrow & \tilde{X} \\ \downarrow & & {\scriptstyle p}\downarrow \\ \mathcal{D} & \longrightarrow & X \end{array}$$

is a pushout diagram (and the same is true after restricting to an open subset of $X$). Exactly as in the proof of the semi-ample part of (2.6), when we replace sections by powers, $\sigma^D$ extends to $\sigma^{\mathcal{D}}$, and $\sigma^{\tilde{X}}$ and $\sigma^{\mathcal{D}}$ have the same pullbacks to $\mathcal{C}$. Thus by (2.11) and (2.7), the sections join to give a section of $L$, nonvanishing at $x$. □

2.12 COROLLARY. *Let $X$ be a scheme, projective over a field of positive characteristic (resp. a finite field). Assume $X$ is a union of closed subsets $X = X_1 \cup X_2$. Let $L$ be a nef line bundle on $X$ such that $L|_{X_i}$ is EWM (resp. semi-ample) for $i = 1, 2$. Let $g_i : X_i \to Z_i$ be the map associated to $L|_{X_i}$. Assume that all but finitely many geometric fibres of $g_2|_{X_1 \cap X_2}$ are connected. Then $L$ is EWM (resp. semi-ample).*

*Proof.* By (2.2) we may assume that the basefield is algebraically closed. Let $G$ be the union of the finitely many fibers of $g_2$ such that $G \cap X_1$ is (nonempty and) not connected. By (2.9), it is enough to show that $X_1 \cup G$ is EWM (resp. semi-ample). Thus we can change notation, and assume from the start that $L|_{X_2}$ is numerically trivial (resp. torsion by (2.16)). But in this case $g_i(X_1 \cap X_2)$ is zero-dimensional, for $i = 1, 2$; thus we can reverse the factors, repeat the argument, and reduce to the case when $L$ is numerically trivial. In this situation the result is obvious (resp., follows from (2.16)). □

2.12.1 *Remark.* Note that the connectedness condition of (2.12) holds in particular if $L$ is numerically trivial on $X_1 \cap X_2$, since then $g_2|_{X_1 \cap X_2}$ has only finitely many geometric fibres.

2.13 COROLLARY. *Let $X$ be a reduced one-dimensional scheme, projective over a field of positive characteristic (resp. a finite field). Any nef line bundle on $X$ is* EWM *(resp., semi-ample).*

*Proof.* Let $L$ be a nef line bundle on $X$. Let $X_1$ be the union of irreducible components on which $L$ is numerically trivial. Let $X_2$ be the union of the remaining components. $L|_{X_2}$ is ample and $L|_{X_1}$ is EWM (resp. torsion by (2.16)). Now apply (2.12). □

2.14 COROLLARY (notation as in (2.3)). *Assume the basefield is perfect, and the spaces are projective. Let $D \subset X$ be a reduced subscheme and let $C \subset Y$ be the reduction of its inverse image. Assume $D$ contains $\mathcal{D}$ set-theoretically. Assume $L|_Y$ and $L|_D$ are* EWM *(resp. semi-ample, and the basefield is finite). If all but finitely many geometric fibres of $g_Y|_C$ are geometrically connected, then $L$ is* EWM *(resp. semi-ample).*

*Proof.* Let $G \subset Y$ be the union of the finite fibres of $g_Y$ such that $G \cap C$ is (nonempty and) not geometrically connected. Let $C' = C \cup G$. Now $L|_{C'}$ is EWM (resp. semi-ample), by (2.12.1). Clearly $g_Y|_{C'}$ has geometrically connected fibres, so (2.6) applies. □

2.15 COROLLARY. *Let $T$ be a reduced purely two-dimensional scheme, projective over an algebraically closed field of positive characteristic (resp. the algebraic closure of a finite field). Let $L$ be a nef line bundle whose restriction to each irreducible component of $T$ has numerical dimension one. Let $p: \tilde{T} \to T$ be the normalization, and let $C \subset \tilde{T}$ be the reduction of the conductor. Assume $L|_{\tilde{T}}$ is* EWM *(resp. semi-ample). Assume further that, set-theoretically, $C$ meets any generic fibre of the associated map in at most one point. Then $L$ is* EWM *(resp. semi-ample).*

*Proof.* This is immediate from (2.14) and the pushout diagram (2.11). □

2.16 LEMMA. *Let $E$ be a projective scheme over the algebraic closure of a finite field. Any numerically trivial line bundle on $E$ is torsion.*

*Proof.* $E$, and any fixed line bundle, are defined over a finite field $k$. By [AK80] the Picard functor $\operatorname{Pic}_{E/k}^H$, for a fixed Hilbert polynomial $H$, is coarsely represented by an algebraic space, $P^H$, proper (in particular of finite type) over $k$. By the Riemann-Roch theorem, [Ful84,18.3.1], any numerically trivial line bundle has constant Hilbert polynomial $\chi(\mathcal{O}_E)$. Thus, since $P^{\chi(\mathcal{O}_E)}$ has only finitely many $k$-points (it is of finite type), the group of numerically trivial line bundles (defined over $k$) is finite. In particular any such line bundle is torsion. □

## 3. Counterexamples in characteristic zero

*Contracting the diagonal of $C \times C$.* Throughout this section, we work over a basefield $k$.

*Notation for Section* 3. Let $C$ be a curve of genus $g$ at least 2, $S = C \times C$, $\pi_i$, the two projections, $\Delta \subset S$ the diagonal, and $L = \omega_{\pi_1}(\Delta)$. Let $I = I_\Delta$.

3.0 THEOREM. *$L$ is* nef *and* big. *If the characteristic of the basefield is positive, $L$ is semi-ample, but in characteristic zero, $L$ is not semi-ample.* (*In any characteristic*) *$\omega(2\Delta)$ is semi-ample, and defines a birational contraction of $\Delta$ to a projective Gorenstein surface, with ample canonical bundle.*

3.0.1 DEFINITION-LEMMA. *There is a $\mathbb{Q}$-line bundle $\omega_\pi$ on the coarse moduli space $\overline{M}_{g,1}$, such that for any family $f : W \to B$ of stable curves, and the induced map $j : W \to \overline{M}_{g,1}$, $\omega_f$ and $j^*(\omega_\pi)$ agree in $\operatorname{Pic}(W) \otimes \mathbb{Q}$.*

*Proof.* See [Mu77]. □

3.1 COROLLARY. *In characteristic zero, $\omega_\pi \in \operatorname{Pic}(\overline{M}_g) \otimes \mathbb{Q}$ is nef and big, but not semi-ample, for any $g \geq 3$.*

*Proof.* Nefness and bigness are instances of (4.4).
Let $E$ be an elliptic curve. Let $T = C \times E$. Let $W$ be the surface (with local complete intersection singularities) obtained by gluing $S$ to $T$ by identifying $\Delta \subset S$ with the horizontal section $C \times p \subset T$. The first projections on each component induces a family of stable curves of genus $g+1$, $f : W \to C$ and $\omega_f|_S = L$; see (5.3). Thus by (3.0), $\omega_\pi$ cannot be semi-ample. □

3.2 LEMMA. *$L|_\Delta$ is trivial, $c_1(L)^2 > 0$, and $L \cdot D > 0$ for any irreducible curve $D \subset S$ other than $\Delta$. The same hold for $\omega(2\Delta)$.*

*Proof.* These are easy calculations using the adjunction formula. □

Let $\Delta_k$ be the $k^{\text{th}}$ order neighborhood of $\Delta \subset S$. Note that the map

$$\omega_C \xrightarrow{\pi_1^* - \pi_2^*} \Omega_S^1 \otimes \mathcal{O}_\Delta$$

induces an isomorphism $j : \omega_C \to I/I^2$.

3.3 LEMMA.
(1) *There is an exact sequence*

$$0 \longrightarrow H^1(I/I^2) \longrightarrow \operatorname{Pic}(\Delta_2) \longrightarrow \operatorname{Pic}(\Delta) \longrightarrow 0.$$

(2) *Let $h : \operatorname{Pic}(C) \to H^{1,1}(C)$ be the map*

$$M \to [\pi_1^*(M) \otimes \pi_2^*(M^*)] \in H^1(I/I^2).$$

*Let*

$$g = j^{-1} \circ h : \operatorname{Pic}(C) \to H^{1,1}(C),$$

$g(M) = c_1(M)$.

(3) *Let $\psi \in \operatorname{Aut}(S)$ switch the two factors. $\psi$ acts on $H^1(I/I^2)$ by multiplication by $-1$. $V \in \operatorname{Pic}(\Delta_2) \otimes \mathbb{Q}$ is fixed by $\psi$ if and only if $V = \pi_1^*(L) \otimes \pi_2^*(L)$ for some $L \in \operatorname{Pic}(C) \otimes \mathbb{Q}$.*

*Proof.* (1) follows from the exact sequence of sheaves of abelian groups

$$1 \to I/I^2 \xrightarrow{x \to 1+x} \mathcal{O}_{\Delta_2}^* \to \mathcal{O}_\Delta^* \to 1;$$

see [H77,III.4.6].

For (2), it is enough to check the result for $M = \mathcal{O}(P)$ for a point $P \in C$. Let $P \in U \subset C$ be an affine neighborhood of $P$, such that $P$ is cut out by $z \in \mathcal{O}_C(U)$. Let $V = C \setminus P$. Let $\mathcal{U}$ be the open cover $\{U \times U, V \times V\} \cap \Delta$, of $\Delta \subset S$. In $\operatorname{Pic}(\Delta_2)$, $\pi_1^*(\mathcal{O}(P)) \otimes \pi_2^*(\mathcal{O}(-P))$ is represented by the cocycle $z_1/z_2 \in H^1(\mathcal{U}, \mathcal{O}_{\Delta_2}^*)$, where $z_i = \pi_i^*(z)$. Under the exact sequence (1) this corresponds to the cocycle

$$\phi = \frac{z_1 - z_2}{z_2} \in H^1(\mathcal{U}, I/I^2).$$

Under the inclusion $I/I^2 \xrightarrow{d} \Omega_S^1 \otimes \mathcal{O}_\Delta$, $\phi$ maps to the cocycle

$$\frac{d(z_1 - z_2)}{z_2} = d(z_1)/z_1 - d(z_2)/z_2$$
$$= (\pi_1^* - \pi_2^*)(d(z)/z)$$

of $H^1(\mathcal{U}, \Omega_S^1|_\Delta)$ (note $z_1$ and $z_2$ are the same in $\mathcal{O}_\Delta$). Thus $g(\mathcal{O}(P))$ is represented by the cocycle $d(z)/z \in H^1(\mathcal{U}, \omega_C)$, which represents $c_1(\mathcal{O}(P))$; see e.g. [H77,III.7.4]. Hence we have (2).

For (3): After tensoring with $\mathbb{Q}$, $\pi_1^* \otimes \pi_2^*$ induces a $\psi$ equivariant splitting of (1). By (2), $\psi$ acts on $H^1(I/I^2)$ by multiplication by $-1$. (3) follows. □

*Remark.* (3.3.2) is a special case of Atiyah's construction of Chern classes (further news from A. Vistoli). (See [At57].)

3.4 LEMMA. *In characteristic zero, $L|_{\Delta_2}$ is nontorsion.*

*Proof.* Suppose $\omega_{\pi_1}(\Delta)|_{\Delta_2}$ is torsion. Then obviously the same is true of $\omega_{\pi_2}(\Delta)$ and thus
$$(\omega_{\pi_2}(\Delta) \otimes \omega_{\pi_1}(\Delta)^*)|_{\Delta_2} = g(\omega_C)$$
is torsion, contradicting (3.3.2). □

3.5 LEMMA. *$\omega(2\Delta)|_{\Delta_k}$ is trivial for all $k \geq 0$.*

*Proof.* There are the usual exact sequences
$$0 \longrightarrow I^k/I^{k+1} \longrightarrow \mathcal{O}_{\Delta_k} \longrightarrow \mathcal{O}_{\Delta_{k-1}} \longrightarrow 0.$$
Also, $H^1(I^k/I^{k+1}) = H^1(\omega_C^{\otimes k})$ is trivial for $k \geq 2$; thus by (3.2) we can assume $k = 1$. Now $\omega(2\Delta)|_{\Delta_1}$ is in $H^1(I/I^2)$ by (3.2), and fixed by $\psi$, thus trivial by (3.3.3). □

*Proof of* (3.0). Assume (in characteristic zero) that $|L^{\otimes n}|$ is basepoint free for some $n > 0$. Then by (3.2), $L^{\otimes n}$ is trivial on $\Delta_2$, contradicting (3.4).

By (3.5) and (3.2), as in the proof of the example theorem, $\omega(2\Delta)$ is semi-ample, and defines a contraction $p : S \to \overline{S}$ of $\Delta$ to a projective normal surface $\overline{S}$. Also, $\omega(2\Delta)$ is the pullback of an ample line bundle, $M$ on $\overline{S}$ and $M = K_{\overline{S}}$, since they agree away from $p(\Delta)$. Thus $\overline{S}$ is Gorenstein, with ample dualizing sheaf. □

3.6 PROPOSITION. *For $g = 2h$, $h \geq 2$, with basefield $\mathbb{C}$, $K_1\mathcal{M}_g^1$ does not carry a scheme-structure so that $f : \overline{M}_{g,1} \to K_1\mathcal{M}_g^1$ is a morphism.*

*Proof.* Fix a curve $C$ of genus $h \geq 2$. Let $p : T \to C$ be the family of stable genus $g = 2h$ curves obtained by gluing together two copies of $\pi_1 : C \times C \to C$ along the diagonal, $\Delta$. Let $S \subset T$ be one of the copies of $C \times C$.

Assume $K_1\mathcal{M}_g^1$ is a scheme, and $f$ is a (necessarily proper) morphism of schemes. The composition
$$j : T \xrightarrow{i} \overline{M}_{g,1} \xrightarrow{f} K_1\mathcal{M}_g^1$$
contracts $\Delta$ and is finite away from $\Delta$. Let $q = j(\Delta)$ have an affine neighborhood. By taking the closure of the inverse image of a general effective principal divisor near $q$, we can find an effective divisor $D \subset \overline{M}_{g,1}$, which restricts to a nontrivial effective Cartier divisor on $i(T)$, disjoint from $i(\Delta)$. Let $D' = i^*(D)$ where $D' \subset T$ is a nontrivial effective Cartier divisor, disjoint from $\Delta$.

We will follow the notation of [AC87] for divisors on $\overline{M}_{g,n}$. Note that $i(T)$ meets only one of the irreducible boundary divisors, $\delta_{h,0}$, whose general element corresponds to a curve with two components, each of genus $h$, one of which is marked. Also, $\delta_{h,0} = \pi^*(\delta_h)$, where $\pi$ is the map $\pi : \overline{M}_{g,1} \to \overline{M}_g$. Thus by the description of $\text{Pic}(\overline{M}_{g,1})$, [AC87,3.1], after possibly replacing $D$ by a positive multiple, the class $D' \in \text{Pic}(T)$ is

$$\omega_p^{\otimes a} \otimes p^*(N)$$

for some integer $a$, and some line bundle $N \in \text{Pic}(C)$. Note that $\omega_p|_S = L$ (of (3.0)). Thus by (3.2), since $D'$ is nontrivial and disjoint from $\Delta$, $a > 0$, and $N$ is trivial. But then $L$ is torsion in a Zariski neighborhood of $\Delta$, contradicting (3.4). □

## 4. The relative dualizing sheaf of the universal stable curve

In this section we work over a basefield $k$. The characteristic is arbitrary except where noted.

*Notation for Section* 4. We will indicate by $\overline{\mathcal{M}}_{g,n}$ the stack of $n$-pointed stable curves, and by $\overline{M}_{g,n}$ its coarse moduli space. Let $\pi_n : \mathcal{U}_{g,n} \to \overline{\mathcal{M}}_{g,n}$ with sections $x_1, \ldots, x_n$ be the universal stable $n$-pointed curve. We will often use the same symbol for a section of a family and the divisor (of the total space) which is its image; we will also use the same symbol for a line bundle on a stack and the associated $\mathbb{Q}$-line bundle on the coarse moduli space (see (3.0.1)). Let $\Sigma \subset \mathcal{U}_{g,n}$ be the union of the $n$ (disjoint) universal sections. Let $L := L_{g,n} := \omega_\pi(\Sigma)$.

The proof of (0.4) goes roughly as follows: First we will show that $L$ is nef and big, and $\mathbb{E}(L) \subset \partial \overline{M}_{g,n+1}$. See (4.9). Then by (0.2) it is enough to show that $L|_{\partial \overline{M}_{g,n+1}}$ is semi-ample. The normalizations of the boundary components are themselves products of various $\overline{M}_{g_i,n_i}$, and $L|_{g,n}$ restricts to *the analogous thing*, (4.4). Thus we can argue by induction.

Dropping the last point $x_{n+1}$ induces the contraction functor $\pi : \overline{\mathcal{M}}_{g,n+1} \to \overline{\mathcal{M}}_{g,n}$, which takes a stable $n+1$-pointed curve $(C, x_1 + \ldots x_{n+1})$ to the stable $n$-pointed curve $(c(C), x_1 + \ldots x_n)$. Here $c : C \to c(C)$ is the stabilization of the $n$-pointed curve $(C, x_1 + \ldots x_n)$. If the latter is stable, $c(C) = C$. Otherwise $c : C \to c(C)$ contracts the irreducible component (necessarily smooth and rational) of $C$ containing $x_{n+1}$ to a point. Also, $\pi$ identifies $\overline{\mathcal{M}}_{g,n+1}$ with $\mathcal{U}_{g,n}$. There is an associated diagram

(4.0)
$$\begin{array}{ccccc}
\mathcal{U}_{g,n+1} & \xrightarrow{c} & \overline{\mathcal{M}}_{g,n+1} \times_{\overline{\mathcal{M}}_{g,n}} \mathcal{U}_{g,n} & \xrightarrow{p_2} & \mathcal{U}_{g,n} \\
\pi_{n+1} \downarrow & & p_1 \downarrow & & \pi_n \downarrow \\
\overline{\mathcal{M}}_{g,n+1} & = & \overline{\mathcal{M}}_{g,n+1} & \xrightarrow{\pi} & \overline{\mathcal{M}}_{g,n}
\end{array}$$

where $c$ is the universal contraction. We will use this identification of $\overline{\mathcal{M}}_{g,n+1}$ with $\mathcal{U}_{g,n}$ repeatedly. For further details, see [Ke92, p. 547].

*Stratification by topological type.* We will use the orbifold stratification of $\overline{M}_{g,n}$ by topological type, which we now recall, following (with some adjustments of notation) [HL96, §4].

Let $B_i \subset \overline{M}_{g,m}$ be the locus of curves with at least $i$ singular points. $B_i \subset \overline{M}_{g,m}$ has pure codimension $i$. The stratum $B_i^0 = B_i \setminus B_{i+1}$ parametrizes curves with exactly $i$ singular points. The connected components of $B_i^0$ correspond to topological types, i.e. equivalence classes under equisingular deformation. $B_i^0$ is an orbifold (all we will use is that it is normal).

There is a finite branched cover, $g : \tilde{B}_i \to B_i$, a product of various $\overline{M}_{i,j}$, which we will now describe. The map $g : \tilde{B}_i \to B_i$ is given roughly by normalizing a curve. Understanding $g^*(L|_{B_i})$ is the main point in the proof of (0.4) and we will go into some detail.

Let $T$ be a curve of fixed topological type, corresponding to a connected component of $B_i^0$. Let $p : \tilde{T} \to T$ be the normalization of $T$. The labeled points $x_i$ are smooth points of $T$, and so give points of $\tilde{T}$. Let $T_1, \ldots, T_v$ be the irreducible components of $T$, with normalizations $\tilde{T}_1, \tilde{T}_2, \ldots \tilde{T}_v$. Let $Y \subset \tilde{T}$ be $p^{-1}(\mathrm{Sing}(T))$. Choose some ordering of the points of $Y$. There is a fixed point-free involution, $\sigma$, of $Y$ so that $T$ is recovered from $\tilde{T}$ by identifying points of $Y$ according to the involution. Let $X \subset T$ be the union of the $x_i$. Let $g_i$ be the genus of $\tilde{T}_i$, and $X_i, Y_i$ the intersections of $X, Y$ with $\tilde{T}_i$. Let $\overline{\mathcal{M}}_T^{X \cup Y}$ be the product

$$\overline{\mathcal{M}}_T^{X \cup Y} := \underset{1 \leq i \leq v}{\times} \overline{\mathcal{M}}_{g_i, X_i \cup Y_i}$$

and $p_i$ the projection onto the $i^{\text{th}}$ component. Here $\overline{\mathcal{M}}_{g_i, X_i \cup Y_i} := \overline{\mathcal{M}}_{g_i, r_i}$ where $r_i$ is the cardinality of $X_i \cup Y_i$ (this alternative notation has the advantage of indicating how the points are distributed). Let $\overline{M}_T^{X \cup Y}$ be the associated coarse moduli space. Let

$$\mathcal{U}_i := p_i^*(\mathcal{U}_{g_i, X_i \cup Y_i})$$

be the pullback of the universal family

$$\pi : \mathcal{U}_{g_i, X_i \cup Y_i} \to \overline{\mathcal{M}}_{g_i, X_i \cup Y_i}$$

from the $i^{\text{th}}$ component of $\overline{\mathcal{M}}_T^{X \cup Y}$. Construct a family $\mathcal{F}_T \to \overline{\mathcal{M}}_T^{X \cup Y}$ by gluing together the $\mathcal{U}_i$ along sections $y_j$ according to the involution $\sigma$. There is a quotient map

$$q : \coprod \mathcal{U}_i \to \mathcal{F}_T,$$

where $\mathcal{F}_T \to \overline{\mathcal{M}}_T^{X \cup Y}$ is a family of stable $m$-pointed curves, and so induces a map $g : \overline{\mathcal{M}}_T^{X \cup Y} \to \overline{\mathcal{M}}_{g,m}$ and a commutative diagram

(4.1)
$$\begin{array}{ccccc} \coprod \mathcal{U}_i & \xrightarrow{q} & \mathcal{F}_T & \longrightarrow & \mathcal{U}_{g,m} \\ \downarrow & & \downarrow & & \downarrow \pi \\ \overline{\mathcal{M}}_T^{X \cup Y} & = & \overline{\mathcal{M}}_T^{X \cup Y} & \xrightarrow{g} & \overline{\mathcal{M}}_{g,m} \end{array}$$

where the right-hand square is fibral.

Note that $\mathcal{F}_T$ has ordinary double points, $q : \coprod \mathcal{U}_i \to \mathcal{F}_T$ is the normalization, and the restriction of the conductor to $\mathcal{U}_i$ is (the union of sections) $Y_i$. Thus by (5.3), if $h : \mathcal{U}_i \to \mathcal{U}_{g,m}$ is given by the top row of (4.1), then

(4.2) $$h^* \omega_\pi(x_1 + \ldots x_m) = p_i^* \omega_\pi(X_i + Y_i).$$

Let $\tilde{B}_i$ be the disjoint union of the $\overline{M}_T^{X \cup Y}$, over the possible topological types $T$ (with $i$ singular points). Now $g$ induces a finite surjective map $g : \tilde{B}_i \to B_i$. Let $\partial \overline{M}_T^{X \cup Y} \subset \overline{M}_T^{X \cup Y}$ be the union of the inverse images of the boundaries $\partial \overline{M}_{g_i, X_i \cup Y_i}$ under the projections onto each component. Let $\partial \tilde{B}_i \subset \tilde{B}_i$ be the union of the $\partial \overline{M}_T^{X \cup Y}$. Set-theoretically,

(4.3) $$\partial \tilde{B}_i = g^{-1}(B_{i+1}).$$

The connected component of $B_i \setminus B_{i+1}$ corresponding to $T$ coarsely represents the stack obtained from $\overline{\mathcal{M}}_T^{X \cup Y} \setminus \partial \overline{\mathcal{M}}_T^{X \cup Y}$ by forgetting the ordering on the irreducible components. This stack is the quotient stack for the action of a finite group (some product of symmetric groups) on $\overline{\mathcal{M}}_T^{X \cup Y}$, and so in particular is smooth. Thus $B_i \setminus B_{i+1}$ is an orbifold. All we will use is:

4.3.1 LEMMA. $B_i \setminus B_{i+1}$ is normal.

Now consider a topological type $T$ of an $n+1$-pointed curve. Suppose $x_{n+1} \in T_j$. Let $c : T \to c(T)$ be the contraction obtained by dropping $x_{n+1}$. There is a corresponding map

$$p_n : \overline{\mathcal{M}}_T^{X \cup Y} \to \overline{\mathcal{M}}_{c(T)}^{X \cup Y \setminus x_{n+1}}$$

and a commutative diagram

$$\begin{array}{ccc} \overline{\mathcal{M}}_T^{X \cup Y} & \xrightarrow{g} & \overline{\mathcal{M}}_{g,n+1} \\ p_n \downarrow & & \downarrow \pi \\ \overline{\mathcal{M}}_{c(T)}^{X \cup Y \setminus x_{n+1}} & \xrightarrow{g} & \overline{\mathcal{M}}_{g,n}. \end{array}$$

All but the $j^{\text{th}}$ components of $\overline{\mathcal{M}}_T^{X \cup Y}$ and $\overline{\mathcal{M}}_{c(T)}^{X \cup Y \setminus x_{n+1}}$ are the same.

If $(T, x_1, \ldots, x_n)$ is stable, then $c : T \to c(T)$ is the identity on underlying curves, and $p_n$ identifies $\overline{\mathcal{M}}_T^{X \cup Y}$ with $\mathcal{U}_j$, and $g : \overline{\mathcal{M}}_T^{X \cup Y} \to \overline{\mathcal{M}}_{g,n+1}$ with $h$ of (4.2).

If $(T, x_1, \ldots, x_n)$ is unstable, then the $j^{\text{th}}$ component of $\overline{\mathcal{M}}_T^{X \cup Y}$ is a point, $\overline{\mathcal{M}}_{0,3}$, and $p_n$ is an isomorphism. Let $p = c(x_{n+1}) \in c(T)$. Either $p$ is a singular point, or one of the labeled points. In the first case choose some irreducible component $p \in c(T)_s$ and a point $y_i$ in the normalization mapping to $p$. Then $p_n$ identifies $\overline{\mathcal{M}}_T^{X \cup Y}$ with the section $y_i$ of $\mathcal{U}_s$. In the second case, say $p = x_i \in c(T)_s$. Then $p_n$ identifies $\overline{\mathcal{M}}_T^{X \cup Y}$ with the section $x_i$ of $\mathcal{U}_s$, and $g$ with the restriction of $h$.

4.4 LEMMA (notation as above).

$$L_{g,n}|_{\overline{M}_T^{X \cup Y}} = p_j^*(L_{g_j, r_j - 1})$$

where $L_{0,2}$ indicates the trivial line bundle on the zero-dimensional space $\overline{\mathcal{M}}_{0,3}$.

*Proof.* This is immediate from the above identifications and (4.2), when the adjunction formula is used in the unstable case to show that the left-hand side is trivial. □

4.5 LEMMA. $L_{0,3}$ and $L_{1,1}$ are semi-ample and big.

*Proof.* $\mathcal{U}_{3,0} = \mathbb{P}^1$ and $L_{0,3} = \mathcal{O}_{\mathbb{P}^1}(1)$. The claims for $L_{1,1}$ are easily checked by consideration of the family of pointed elliptic curves given by a general pencil of plane cubics. □

4.6 LEMMA (notation as in (4.0)).

$$\omega_{\pi_{n+1}}(x_1 + \cdots + x_n) = (p_2 \circ c)^*(\omega_{\pi_n}(x_1 + \cdots + x_n)).$$

*Proof.* The left-hand side is $\pi_{n+1}$ nef; in particular, it is $c$ nef.

$$(\overline{\mathcal{M}}_{g,n+1} \times_{\overline{\mathcal{M}}_{g,n}} \mathcal{U}_{g,n}, x_1 + \cdots + x_n)$$

has canonical singularities, by inversion of adjunction ([Kol96b,7.5]) (or explicit local coordinates; see [Kn83]). Also, $c$ is birational, so the result follows by negativity of contractions; see [Kol et al.92,2.19]. □

4.7 PROPOSITION. $L_{g,n}$ *is nef.*

*Proof.* First we induct on $n$ to reduce either to (4.5), or the case $g \geq 2$ and $n = 0$.

Assume $L_{g,n}$ is nef. By (4.6)

(4.7.1) $$L_{g,n+1} = (p_2 \circ c)^*(L_{g,n}) + x_{n+1}.$$

Since $L_{g,n+1}|_{x_{n+1}}$ is trivial by adjunction, $L_{g,n+1}$ is nef. Thus we may assume that $g \geq 2$ and $n = 0$.

4.7.2 *Remark* (to be used in the proof of (4.8)).

(4.7.3) $$p_2 \circ c \circ x_{n+1} : \overline{\mathcal{M}}_{g,n+1} \to \mathcal{U}_{g,n}$$

is an isomorphism; thus if $d = 3g - 1 + n$ is the dimension of $\mathcal{U}_{g,n+1}$, then

$$c_1(L_{g,n+1})^d \geq c_1((p_2 \circ c)^*(L_{g,n}))^{d-1} \cdot x_{n+1} = c_1(L_{g,n})^{d-1}.$$

Hence if $L_{g,n}$ is big, then so is $L_{g,n+1}$.

For (4.7) it is enough to consider a one-dimensional family of stable curves. By (4.4) and the above reduction, we can assume $g \geq 2$, $n = 0$ and the general fibre is smooth. Now the result follows from [Kol90,4.6]. □

4.8 PROPOSITION. *Let $p : \mathcal{C} \to B$ with sections $x_1, \ldots, x_n$ be a family of stable n-pointed curves of genus g. Assume B is irreducible, the general fibre of p is smooth, and the associated map $B \to \overline{M}_{g,n}$ is generically finite. If $V \subset \mathcal{C}$ is a subvariety surjecting onto B, and V is not contained in any of the $x_i$, then $\omega_p(x_1 + \ldots x_n)|_V$ is big.*

*Proof.* Note that there are two cases. Either $V = \mathcal{C}$ or $p|_V$ is generically finite.

Consider the first case. We start by inducting on $n$, to reduce to the case $g \geq 2$, and $n = 0$, as in the proof of (4.7). Suppose there are $n + 1$ sections (and the result is known for an $n$-pointed curve). Let $c : \mathcal{C} \to \mathcal{E}$ be the contraction of $x_{n+1}$. This induces a map $B \to \overline{M}_{g,n}$, which is just the composition $B \to \overline{M}_{g,n+1} \xrightarrow{\pi} \overline{M}_{g,n}$. After replacing $B$ by a finite cover (to deal with the fact that $\overline{M}_{g,n}$ does not carry a universal family, see [V95,9.25]), we can assume there is a map $B \to B'$ so that $\mathcal{E} \to B$ is pulled back from an $n$-pointed curve $\mathcal{D} \to B'$, for which the induced map $B' \to \overline{M}_{g,n}$ is generically finite. We have a commutative diagram

$$\begin{array}{ccccccc}
\mathcal{C} & \xrightarrow{c} & \mathcal{E} = B \times_{B'} \mathcal{D} & \xrightarrow{p_2} & \mathcal{D} & \longrightarrow & \mathcal{U}_{g,n} \\
\downarrow & & p_1 \downarrow & & \downarrow & & \pi_n \downarrow \\
B & = & B & \longrightarrow & B' & \longrightarrow & \overline{\mathcal{M}}_{g,n}
\end{array}$$

where the right two squares are fibral, and the left-most square is pulled back from the left-most square of (4.0). There are formulae analogous to (4.7.1), (4.7.3), so we can apply induction exactly as in (4.7.2).

We may assume $g \geq 2$ and $n = 0$. By [LT89,4.8] (or [V77,2.10] in characteristic zero)

$$\omega_p = p^* c_1(p_*(\omega_p)) + Z$$

where $Z$ is an effective Cartier divisor, surjecting onto $B$ (the intersection of $Z$ with the generic fibre is the set of Weierstrass points). Also, $c_1(p_*(\omega_p))$ is big by [C93,2.2]. Thus by Kodaira's lemma, we may write $c_1(p_*(\omega_p)) = A + E$ with $A$ ample, and $E$ effective and Cartier. Let $d-1$ be the dimension of $B$. Note that $\omega_p$ is nef by (4.7).

$$c_1(\omega_p)^d \geq p^*(A)^{d-1} \cdot Z > 0.$$

Thus $\omega_p$ is big.

Now consider the second case. After pulling back, we can assume $V = \sigma(B)$ for a section $\sigma$, distinct from the $x_i$. We can also assume that $B$ is normal. Thus $\mathcal{C}$ is normal.

By (4.7) and the first case, $\omega_p(\Sigma)$ is nef and big. Thus by Kodaira's lemma

$$\omega_p(\Sigma) = A + E$$

with $A$ ample and $E$ effective and Cartier.

Suppose $\omega_p(\Sigma)|_\sigma$ is not big. Then by (1.2), $\sigma$ is in the support of $E$. Thus there exists $\lambda > 0$ so that $\lambda E = \sigma + V$, with $V$ an effective $\mathbb{Q}$-Weil divisor, whose support does not contain $\sigma$. Then

$$(1+\lambda)\omega_p(\Sigma)|_\sigma = \omega_p(\sigma + V)|_\sigma + (\lambda A + \Sigma)|_\sigma$$

and the first term (on the right-hand side of the above equality) is effective by adjunction; see (5.3). The second is big. Thus $\omega_p(\Sigma)|_\sigma$ is big, which is a contradiction. □

4.9 COROLLARY. *$L_{g,n}$ is nef and big, and its exceptional locus is contained in $\partial \overline{M}_{g,n+1}$.*

4.10 THEOREM. *If the basefield has positive characteristic, then $L$ is semi-ample.*

*Proof.* We will proceed by induction on $m = 3g - 2 + n$, the dimension of $\mathcal{U}_{g,n}$. Note that $\mathcal{U}_{g,n}$ and $L$ are defined over the characteristic field; thus (4.10) holds for $m \leq 2$ by (0.3). By (4.9) and (0.2) it is enough to show $L|_{\partial \overline{M}_{g,n+1}}$ is semi-ample.

We will prove that $L|_{B_i}$ is semi-ample by induction on $i$. Of course $B_i$ is empty, and there is nothing to prove, for $i > m$. By (4.3),(4.4) and (4.9), $\mathbb{E}(L|_{\tilde{B}_i}) \subset \partial \tilde{B}_i$. By (4.9) and induction on $m$, $L|_{\tilde{B}_i}$ is semi-ample. Thus $L|_{B_i}$ is semi-ample by induction on $i$, (4.3.1) and (2.10.1). □

## 5. Existence of birational $K + \Delta$ negative extremal contractions on 3-folds of positive characteristic

5.0 *Proof of* (0.5). By Kodaira's lemma $L = A + E$ for $A$ ample, and $E$ effective. Write $E = N_0 + N_1 + N_2$, where $N_i$ is the sum of the irreducible components (with the same coefficients as in $E$) on which $L$ has numerical dimension $i$. Let $T$ be the support of $N_1$.

The main issue will be to show that $L|_T$ is EWM (resp. semi-ample if the basefield is finite). Suppose this has been established. Let $R_i$ be the support of $N_i$ for $i = 0, 2$. Let $W = R_0 \cup \mathbb{E}(L|_{R_2})$. By (1.2), $\mathbb{E}(L) = T \cup W$. Note that $\mathbb{E}(L|_{R_2})$ is one-dimensional, so $L|_W$ is numerically trivial (resp. torsion, by (2.16)). Thus $L|_{\mathbb{E}(L)}$ is EWM (resp. semi-ample) by (2.12.1), and so $L$ is EWM (resp. semi-ample) by (0.2).

Now we consider $L|_T$. Let $T = \bigcup T_i$ be the decomposition into irreducible components. Let $p : \tilde{T} \to T$ be the normalization, and $\tilde{T} = \bigcup \tilde{T}_i$ the corresponding decomposition into connected components. We will show first that $L|_{\tilde{T}}$ is semi-ample. Since $L|_{\tilde{T}_i}$ has numerical dimension one, it is enough to show $L|_{\tilde{T}_i}$ moves; see (5.2). This will follow from a simple Riemann-Roch calculation. It will follow from adjunction that the conductor of $\tilde{T}_i \to T_i$ is generically a section for the map associated to $L|_{\tilde{T}_i}$ (and something slightly weaker holds for $p : \tilde{T} \to T$). Then (2.12) and (2.15) will imply $L|_T$ is EWM (resp. semi-ample).

By (2.2) we may assume the basefield is algebraically closed.

We start with some bookkeeping, to prepare for adjunction. Let $M = L - (K_X + \Delta)$, which is nef and big by assumption. Let

$$N = \sum_{i=1}^{r} a_i T_i$$

be the decomposition into irreducible components. Define $\lambda_i > 0$ by $\lambda_i a_i + e_i = 1$, where $e_i$ is the coefficient of $T_i$ in $\Delta$. Arrange indices so that $\lambda_1 \geq \lambda_2 \ldots \geq \lambda_r$. Define $\Gamma_i$ by

$$\Delta + \lambda_i E = T_i + \sum_{j > i} T_i + \Gamma_i.$$

By construction, $\Gamma_i$ is effective, its support does not contain $T_i$ and

(5.0.1) $$(1 + \lambda_i)L - (K_X + \sum_{j \geq i} T_j + \Gamma_i) = M + \lambda_i A := A_i$$

is ample.

Let $C_i \subset \tilde{T}_i$ and $D_i \subset T_i$ be the conductors. Let $T^i = \bigcup_{j \geq i} T_j$, with normalization $\tilde{T}^i$ and conductor $C^i \subset \tilde{T}^i$. Let $Q_i \subset \tilde{T}_i$ be the restriction of $C^i$ (note this is just the restriction to a connected component). We will use the same symbol to indicate the integral Weil divisor associated to each conductor. $C_i$ is a subscheme of $Q_i$; thus there is an inequality between Weil divisors $Q_i \geq C_i$.

By the adjunction formula, (5.3), below,

(5.0.2) $$(K_X + \sum_{j \geq i} T_j + \Gamma_i)|_{\tilde{T}_i} = K_{\tilde{T}_i} + Q_i + R_i$$

for some effective $R_i$. Of course (5.0.1) and (5.0.2) imply

(5.0.3) $$K_{\tilde{T}_i} + Q_i + R_i = (1+\lambda_i)L|_{\tilde{T}_i} - A_i|_{\tilde{T}_i}.$$

By the Riemann-Roch theorem, the leading term of $\chi(L^{\otimes r} \otimes \mathcal{O}_{\tilde{T}_i})$ is

$$r/2(L|_{\tilde{T}_i}) \cdot (L|_{\tilde{T}_i} - K_{\tilde{T}_i}).$$

Using (5.0.3) and the fact that $L|_{\tilde{T}_i}$ has numerical dimension one we get

$$\begin{aligned} r/2(L|_{\tilde{T}_i}) \cdot (L|_{\tilde{T}_i} - K_{\tilde{T}_i}) &= r/2(L|_{\tilde{T}_i}) \cdot (-K_{\tilde{T}_i}) \\ &= r/2(L|_{\tilde{T}_i}) \cdot (Q_i + R_i + A_i|_{\tilde{T}_i}) \\ &\geq r/2(L \cdot A_i \cdot T_i). \end{aligned}$$

$L \cdot A_i \cdot T_i$ is strictly positive by the Hodge index theorem. Thus $L|_{\tilde{T}}$ is semi-ample by (5.4) and (5.2). Let $f_i$ be the map associated to $L|_{\tilde{T}_i}$.

By (5.0.3), $K_{\tilde{T}_i} + Q_i$ is negative on any generic fibre of $f_i$. Now, $Q_i$ and $C_i$ have integral coefficients. Thus, by adjunction, any generic fibre is $\mathbb{P}^1$, and in a neighborhood of any generic fibre, either $Q_i$ is empty, or a section. Since $Q_i \geq C_i$, the same holds for $C_i$.

Thus $L|_{T_i}$ is EWM (resp. semi-ample) by (2.15). Let $g_i : T_i \to Z_i$ be the associated map.

Let

$$W' = \bigcup_{i=1}^{i=n-1} T_i,$$

with reduced structure. We will argue inductively that $L|_{W'}$ is EWM (resp. semi-ample). Assume this holds for $n-1$, with associated map $g : W' \to Z$, and consider $T_n \cup W'$. We claim that $T_n \cap W'$ meets any generic fibre of $g$ in at most one (set-theoretic) point. The theorem follows from the claim, by (2.12).

To see the claim, suppose on the contrary, that there are points $p \neq q$ of $T_n$ along some generic fibre, $G$. Let $X \subset G$ be a minimal connected union of irreducible components, containing $p, q$ (minimal under inclusion). Let $F$ be an irreducible component of $X$ which lies on a $T_i$ with $i$ minimal. By minimality of $X$,

$$F \cap (\{p,q\} \cup \text{Sing}(X))$$

contains at least two points. Thus by the minimality of $i$, there are at least two distinct singular points of $T^i$ along $F$, and thus $Q_i$ meets the strict transform $\tilde{F} \subset \tilde{T}_i$ in at least two distinct points. Observe that $\tilde{F}$ is a generic fibre of $f_i$. But by (5.0.3),

$$(K_{\tilde{T}_i} + Q_i)|_{\tilde{F}} = K_{\tilde{F}} + Q_i|_{\tilde{F}} = K_{\mathbb{P}^1} + Q_i|_{\tilde{F}}$$

is negative, which is a contradiction. □

Now we turn to the lemmas used in the proof of (0.5):

5.2 LEMMA. *Let $T$ be a normal surface, projective over an algebraically closed field. Let $L$ be a nef line bundle on $T$, of numerical dimension one. If $h^0(L^{\otimes m}) > 0$ for some $m > 0$, then $L$ is semi-ample.*

*Proof.* Passing to a desingularization, we can assume $T$ is nonsingular. In this case the result is familiar; see, for example, the proof of [Kol et al.92, 11.3.1]. □

5.2.1 *Example.* If one assumes only that $T$ is integral, (5.2) fails. K. Matsuki showed me the following counterexample: Let $C$ be a curve, and let $T$ be obtained from $C \times C$ by gluing together two points on different fibres $F_1, F_2$ of the first projection. Take $L = \mathcal{O}(2F_1 + F_2)$ on $T$.

5.3 ADJUNCTION FORMULA. *Let $T \subset X$ be a reduced Weil divisor on a normal variety $X$. Let $\pi : \tilde{T} \to T$ be the normalization, and let $p$ be the composition $p : \tilde{T} \to T \subset X$. Suppose $K_X + T$ is $\mathbb{Q}$-Cartier. Let $C \subset \tilde{T}$ be the Weil divisor defined by the conductor. Then there is a canonically defined effective $\mathbb{Q}$-Weil divisor $\mathcal{D}$ on $\tilde{T}$, whose support is contained in $p^{-1}(\mathrm{Sing}(X))$, such that*

$$K_{\tilde{T}} + C + \mathcal{D} = p^*(K_X + T).$$

*Proof.* Since we are working with Weil divisors, we can remove codimension two subsets from $T$ and so may assume $T$ is Cohen-Macaulay, and $\tilde{T}$ is nonsingular. By [R94, p. 17], there is a canonical surjective map $\pi^*(\omega_T) \to \omega_{\tilde{T}}(C)$, which is an isomorphism wherever $T$ is Gorenstein. Also, $\omega_X(T) \otimes \mathcal{O}_T = \omega_T$. Thus for each $r$ there is an induced map

$$\omega_{\tilde{T}}(C)^{\otimes r} = p^*(\omega_X(T)^{\otimes r})/\mathrm{torsion} \to p^*(\omega_X(T)^{[r]})/\mathrm{torsion}$$

which is an isomorphism wherever $X$ is nonsingular. The result follows. □

Shokurov gives a formula analogous to (5.3), but without isolating the conductor. (See [Sho91,3.1].) When $T$ is Gorenstein in codimension one, Corti gave the same formula on $T$, with the same proof; see [Kol et al.92, 16.5].

5.3.2 *Question.* We wonder if a form of (5.3) holds on $T$ itself. The analog of $K$ on a reduced scheme $T$ of pure dimension $d$ is $-\frac{1}{2}\tau_{d-1}(\mathcal{O}_T)$, where $\tau$ is the Todd class; see [Ful84,ch. 18]. In the context of (5.3), is

$$(K_X + T)|_T + \frac{1}{2}\tau_{d-1}(\mathcal{O}_T)$$

an effective class? This would be useful for Riemann-Roch calculations, for example, in the proof of (0.5). When $T$ is Gorenstein in codimension one, this holds, and follows from (5.3). But it does not follow in general from (5.3); Serre's inequality on conductor lengths

$$n \leq 2\delta$$

goes in the wrong direction. (See [R94,3.2].)

We note one immediate corollary of (5.3):

5.3.3 COROLLARY. *Let $C \subset S$ be an integral curve on a normal $\mathbb{Q}$-factorial surface, projective over a field. If $(K_S + C) \cdot C < 0$, then $C$ is a smooth curve of genus* 0.

*Proof.* Over any singular point, the conductor has degree at least two. Thus by degree considerations the conductor is empty and $\tilde{C} = C$ has genus 0. □

5.4 LEMMA. *Let $X$ be a projective pure dimensional scheme over a field. Let $L$ be a* nef *line bundle on $X$. Then $H^{\dim(X)}(L^{\otimes r})$ is bounded over $r > 0$.*

*Proof.*
$$H^{\dim(X)}(L^{\otimes r}) = H^0(\omega_X \otimes L^{\otimes -r})^*,$$

where $\omega_X$ is the pre-dualizing sheaf (see [H77,III.7.3]). Thus it is enough to show that for any coherent sheaf $\mathcal{F}$, $H^0(\mathcal{F} \otimes L^{-\otimes r})$ is bounded. For this we will use Grothendieck's Dévissage, [Kol96,VI.2.2], and consider the class of coherent sheaves $\mathcal{F}$ with this property. Follow Kollár's notation and suppose $\mathcal{F} = \mathcal{O}_Z$, for an integral subscheme $Z \subset X$. Now, $H^0(Z, L^{\otimes -r})$ vanishes unless $L^{\otimes r}$ is trivial. The other two conditions of [Kol96,VI.2.1] are immediate. □

*Remark.* For much stronger results than (5.4), see [Fuj82].

5.5 *Cone of Curves.* We will follow the notation of [Kol96,II.4] for notions related to cones. In particular $N^1(X)$ indicates the dual of the Néron-Severi group, with $\mathbb{R}$ coefficients. Also an extremal ray, $R$, of a closed convex cone, is a one dimensional subcone which is extremal; i.e. if $x_1 + x_2 \in R$ then $x_1, x_2 \in R$.

5.5.1 *Definition.* A class $h \in N^1(X)$ is called *ample* (resp. *nef*) if $h$ is a strictly positive (resp. nonnegative) function on $\overline{NE}_1(X) \setminus \{0\}$.

For a line bundle on a projective variety, (5.5.1) agrees with the usual notion of ample, by Kleiman's criterion.

5.5.2 PROPOSITION. *Let $X$ be a $\mathbb{Q}$-factorial normal 3-fold projective over a field (of arbitrary characteristic). Let $\Delta$ be an $\mathbb{R}$-boundary. Let $h \in N^1(X)$ be an ample class such that*

$$\eta := K_X + \Delta + h$$

*is* nef, *and numerically equivalent to an effective* $\mathbb{R}$-*divisor,* $\Gamma$. *Then the extremal subcone supported by* $\eta$,

$$\eta^{\perp} \cap \overline{NE}_1(X),$$

*is generated by the classes of a finite number of curves. Furthermore, suppose*

$$\Gamma = \gamma + \alpha$$

*where* $\alpha$ *is ample, and* $\gamma$ *is effective. Let* $S$ *be the support of* $\Delta + \gamma$. *Then any* $\eta$ *trivial extremal ray is generated either by a curve in* $\mathrm{Sing}(S) \cup \mathrm{Sing}(X)$, *or by a rational curve,* $C$, *satisfying*

$$0 < -(K_X + \Delta) \cdot C \leq 3.$$

*Proof.* Let $R$ be an extremal ray with $\eta \cdot R = 0$. We will argue that $R$ is of the form described, necessarily $\gamma \cdot R < 0$. It follows that $R$ is in the image of

$$\overline{NE}_1(T) \to \overline{NE}_1(X)$$

for some irreducible component $T$ of the support of $\gamma$, with $T \cdot R < 0$ (see the proof of [Kol96,II.4.12]). Let $p : \tilde{T} \to T$ be the minimal desingularization of the normalization of $T$, and let $\pi : \tilde{T} \to X$ be the induced map. Necessarily, $R$ is generated by the image of some $\pi^*(\eta)$ trivial extremal ray, $J$, of $\overline{NE}_1(\tilde{T})$.

Define $a \geq 0$ so that

$$\Delta + a \cdot \gamma = T + E$$

for $E$ an effective $\mathbb{R}$-divisor whose support does not contain $T$. Note that $E|_T$ has support contained in $\mathrm{Sing}(S)$. Let $g = h + a\alpha$. By (5.3)

$$\pi^*(K_X + T) = K_{\tilde{T}} + Q$$

for an effective class $Q$, whose image in $X$ has support contained in $\mathrm{Sing}(S) \cup \mathrm{Sing}(X)$. There is a numerical equality:

$$(1 + a)\pi^*(\eta) = K_{\tilde{T}} + Q + p^*(E|_T) + \pi^*(g).$$

Since $g$ is ample,

$$J \cdot (K_{\tilde{T}} + Q + p^*(E|_T)) < 0.$$

Thus either $J$ is $K_{\tilde{T}}$ negative, or is generated by a curve mapping into $\mathrm{Sing}(S) \bigcup \mathrm{Sing}(X)$. Suppose we are in the first case, but not the second. By Mori's cone theorem (for smooth surfaces), $J$ is generated by a smooth rational curve, $C$, with $0 < -K_{\tilde{T}} \cdot C \leq 3$, and there are finitely many possible $C$'s (though the number will in general depend on $\eta$). Now, $T \cdot C < 0$, and $W \cdot C > 0$ for any other irreducible component, $W$, of $S$; thus

$$0 < -(K_X + \Delta) \cdot C \leq -(K_X + T) \cdot C = -(K_{\tilde{T}} + Q) \cdot C \leq -K_{\tilde{T}} \cdot C \leq 3.$$

5.5.3 *Remark. Kodaira's lemma for $\mathbb{R}$-classes.* If one could show that a nef class $\eta$ with $\eta^{\dim X} > 0$ satisfied Kodaira's lemma, i.e. has an expression $h + E$ with $h$ ample, and $E$ effective, then (5.5.2) would imply that the supported extremal subcone is finite rational polyhedral. In dimension two this is easy; see [Kol96,4.12]. It is also known, in all dimensions, that in characteristic zero the proof, due to Shokurov, is an application of Kawamata-Viehweg vanishing; see [Sho96]. This generalized Kodaira lemma has some other interesting implications; see [KMcK96, §2].

*Proof of* (0.6). Let $K_X + \Delta = \gamma$, for an effective class $\gamma$. Let $R$ be an extremal ray, with $(K_X + \Delta) \cdot R < 0$. Let $\eta \in N^1(X)$ be a nef class supporting $R$. Then by compactness of a slice of $\overline{NE}_1(X)$ ([Kol96,II.4.8]), after possibly replacing $\eta$ by a positive multiple,

$$h := \eta - (K_X + \Delta)$$

is ample. (1) through (3) now follow from (5.5.2); see, for example, the proof of [Kol96,III.1.2]. □

5.5.4 *Remark.* It is natural to expect, under some singularity assumptions, that the extremal rays in (0.6) are all generated by smooth rational curves. In characteristic zero this follows from Kawamata-Viehweg vanishing applied to the corresponding extremal contraction. In characteristic $p$, the proof of (0.5), plus some easy analysis in the case when the exceptional locus is a surface contracted to a point, shows that the ray is generated by a (possibly singular) rational curve, except possibly in the case of a small contraction, where, if the singularities are isolated LCIQ and the ray is $K_X$ negative, then any generating curve is rational by [Kol92,6.3].

UNIVERSITY OF TEXAS, AUSTIN, TX
*E-mail address*: keel@math.utexas.edu